\input amstex
\documentstyle{amsppt}
\NoBlackBoxes
%
\chardef\tempcat=\the\catcode`\@
\catcode`\@=11

\font\fourteenbf=cmbx10 scaled \magstep2

\def\QED{\nobreak\quad\ifmmode\roman{Q.E.D.}\else{\rm Q.E.D.}\fi}
\def\keyedby#1{}
\def\affil#1\endaffil{}
\def\ratitle{\setbox0=\hbox{\fourteenbf R}
                 \advance\baselineskip by \ht0  
                 \centerline{\fourteenbf RESEARCH ANNOUNCEMENTS}%
                 \baselineskip45pt}     
\def\date#1\enddate{\thanks Received by the editors #1 \endthanks}
\def\shorttitle#1{\rightheadtext{#1}}
\def\shortauthor#1{\leftheadtext{#1}}

\def\ml #1\endml{\email #1\endemail}
\long\def\ext #1\endext{\block #1\endblock}
\def\ac #1\endac{{\monograph@true \indenti=0pc 
                 \heading Contents\endheading \toc #1\endtoc}}
\def\ch #1\endch{\def\subheadfont@{\it} \subheading{#1} \def\subheadfont@{\bf}}

\let\rem=\remark
\let\endrem=\endremark

\def\RM#1{{\rm#1}}
\def\<{\nobreak\kern-\mathsurround\nobreak}
\everydisplay={\def\.{\thinspace.}}
\def\fighere{\def\figure@type{3}\ins@}

\loadbold 
\def\boldcdot{\boldsymbol\cdot}
\let\varroster=\roster

\def\cprime{\/{\m@th$'$}}

{\let\enddocument=\relax
\gdef\bookrev{%
	\long\gdef\revtop ##1\endrevtop{\def\revtop@{{\interlinepenalty10000
		\tenpoint
		\everypar{\hangindent=\parindent}
                \noindent ##1\endgraf}}}%

	\gdef\reviewer ##1\endreviewer{\def\revr@{{\leftskip0pt plus1fil
	  \rightskip0pt \parfillskip0pt \tenpoint\smc ##1\endgraf}}}%

	\gdef\affil ##1\endaffil{\def\affil@{{\leftskip0pt plus1fil
	  \rightskip0pt \parfillskip0pt \tenpoint\smc ##1\endgraf}}}%

      \gdef\bookrevheading{\hrule height\z@
                 \vskip-\topskip
                 \setbox0=\hbox{\fourteenbf B}
                 \advance\baselineskip by \ht0  
                 \centerline{\fourteenbf BOOK REVIEW}%
                 \baselineskip45pt\hbox{}\baselineskip0pt}

   \gdef\document{%
        \bookrevheading
        \gdef\bookrevheading{\penalty-100\vskip2pc plus2pc}
   	\leftheadtext{BOOK REVIEW}%
	\rightheadtext{BOOK REVIEW}%
	\nobreak
	\vskip8pt  plus5pt minus2pt
	\nobreak
	\revtop@
	\nobreak
	\vskip12pt plus5pt minus2pt
	\tenpoint }%

	\outer\def\enddocument{\nobreak
	\vskip6pt minus6pt
	\nobreak
	\revr@
        \nobreak
        \affil@
	\end }%

}
}

\def\prelogo{}


\def\jourlogo{\hbox{\vbox to0pt{\prelogo
   \vbox{\sixrm \baselineskip 6.5pt
   \parindent 0pt APPEARED IN BULLETIN OF THE\hfil\break
   AMERICAN MATHEMATICAL SOCIETY\hfil\break
   Volume \cvol@, Number \cvolno@, \cmonth@\ \cvolyear@, Pages \cpgs@\endgraf}%
   \vss}}}


\def\cvol#1{\gdef\cvol@{\ignorespaces#1\unskip}}
\def\cvolno#1{\gdef\cvolno@{\ignorespaces#1\unskip}}
\def\cmonth#1{\gdef\cmonth@{\ignorespaces#1\unskip}}
\def\cvolyear#1{\gdef\cvolyear@{\ignorespaces#1\unskip}}
\def\cyear#1{\gdef\cyear@{\ignorespaces#1\unskip}\cyear@@#100000\end@}
\def\cpgs#1{\gdef\cpgs@{\ignorespaces#1\unskip}}

\def\cyear@@#1#2#3#4#5\end@{\gdef\cyearmodc@{#3#4}%
        \gdef\cyearmodcHold@{#3#4}}

\cvol{000}
\cvolno{0}
\cmonth{}
\cyear{0000}
\cvolyear{0000}
\cpgs{}

\font\sixsy=cmsy6

\def\copyrightline@{\baselineskip1.75pc
    \rightline{%
        \vbox{\sixrm \textfont2=\sixsy \baselineskip 7pt
            \halign{\hfil##\cr
                \copyright\cyear@\ American Mathematical Society\cr
                 0273-0979/\cyearmodc@\ \$1.00 + \$.25 per page\cr }}}}

\def\cyearmodc#1{\gdef\cyearmodc@{\ignorespaces#1\unskip}}


\let\logo@=\copyrightline@
\def\pretitle{\jourlogo \vskip24\p@ plus12\p@ minus12\p@}



\def\shoveright#1{\omit\span\omit\span\omit 
  \hfil$\@lign\displaystyle{{}#1}\m@th$%
  \iftagsleft@ 
    \ifx\undefined\displaywidth@ 
    \else\kern-\displaywidth@\fi
  \fi}


\def\aligned@{\bgroup\vspace@\Let@
 \def\shoveright##1{\omit\span\egroup\hfill\bgroup##1}%
 \ifinany@\else\openup\jot\fi\ialign
 \bgroup\hfil\strut@$\m@th\displaystyle{##}$&
 $\m@th\displaystyle{{}##}$\hfil\crcr}


\def\gathered{\null\,\vcenter\bgroup\vspace@\Let@
 \def\shoveright##1{##1\aftergroup\aftergroup\aftergroup\hfilneg}%
 \ifinany@\else\openup\jot\fi\ialign
 \bgroup\hfil\strut@$\m@th\displaystyle{##}$\hfil\crcr}


\def\gather{\RIfMIfI@\DN@{\onlydmatherr@\gather}\else
 \ingather@true\inany@true\def\tag{&}%
 \vspace@\allowdisplaybreak@\displaybreak@\intertext@
 \displ@y\Let@
 \def\shoveright##1{\omit\span\omit
   \hfil\llap{\strut@$\m@th\displaystyle{##1}$}%
   \iftagsleft@
     \ifx\undefined\gdisplaywidth@ 
     \else\kern-\gdisplaywidth@\fi
   \fi}%
 \iftagsleft@\DN@{\csname gather \endcsname}\else
  \DN@{\csname gather \space\endcsname}\fi\fi
 \else\DN@{\onlydmatherr@\gather}\fi\next@}

%
%
\newskip\colspacing      \colspacing=1em
%
%
\def\varquad{\hskip\colspacing\relax}
\def\reduce#1{\advance#1 by -}
%
%
\def\halfspace{\hskip .5\colspacing\relax}
%
%
\def\ruleh{\noalign{\hrule}}
%
%
\def\rulev{\halfspace\vrule\halfspace}%

\def\t@bleside{\ifvmode\hrule \else\vrule\fi\relax}
\newif\iftop@ 
\def\toptablecaption{\top@true\tablecaption}
\def\tablecapf@nt{\tenpoint}
%
%
\def\tablecaption#1{\iftop@\else\medskip\fi\begingroup\tablecapf@nt\Let@
  \centerline{\vbox{\ialign{\hfil##\hfil\cr#1\crcr}}}%
  \endgroup\iftop@\medskip\fi}
%
%
\def\colheading#1{\vbox{\normalbaselines\Let@\ialign{\rm\hfil##\hfil\cr
        \tablestrut depth\z@ #1\unskip\tablestrut height\z@\crcr}}}
%
%
\def\table{\relaxnext@\DN@{\ifx\next\nobox\def\t@bleside{}\fi
 \relaxnext@{}\,\t@bleside\vcenter\bgroup\Let@\vspace@
 \offinterlineskip\t@bleside\ialign\bgroup
  \halfspace\hfil$\tablestrut####$\hfil&&\varquad\hfil$####$\hfil\crcr
}\futurelet\next\next@}
%
\def\nobox{}
\def\endtable{\crcr\egroup\t@bleside\egroup\t@bleside\,}
%

\def\matrix{\null\,\vcenter\bgroup\Let@\vspace@
 \def\ruleh{\noalign{\vskip\tw@\p@\hrule\vskip\thr@@\p@}}%
 \def\everyr@w{}%
 \normalbaselines\openup\spreadmlines@\m@th\ialign
 \bgroup\hfil$##$\hfil&&\varquad\hfil$##$\hfil\crcr
 \Mathstrut@\crcr\noalign{\kern-\baselineskip}}

\def\everyr@w{\halfspace\tablestrut}

\newdimen\rowspacing  


\def\normalbaselines{\lineskip\normallineskip
  \baselineskip\normalbaselineskip \lineskiplimit\normallineskiplimit
  \rowspacing1.6\normalbaselineskip}

%
%
\def\tablestrut{\vrule width\z@ height.67\rowspacing
  depth.33\rowspacing}

\def\deeper#1{\begingroup
  \dimen@.33\rowspacing\advance\dimen@#1\tablestrut depth\dimen@
  \relax\endgroup}
\def\higher#1{\begingroup
  \dimen@.67\rowspacing\advance\dimen@ #1\tablestrut height\dimen@
  \relax\endgroup}

\def\format@#1\\{\def\preamble@{#1}%
 \def\c{\hfil$\the\hashtoks@$\hfil}%
 \def\r{\hfil$\the\hashtoks@$}%
 \def\l{$\the\hashtoks@$\hfil}%
 \def\ctext{\hfil\the\hashtoks@\hfil}%
 \def\rtext{\hfil\the\hashtoks@}%
 \def\ltext{\the\hashtoks@\hfil}%
 \setboxz@h{\xdef\Preamble@{\everyr@w\preamble@}}%
 \ifnum`{=0 \fi\iffalse}\fi
 \ialign\bgroup\span\Preamble@\crcr}

\catcode`\@=\tempcat


\keyedby{bull239/amh}

\topmatter
\cvol{26}
\cvolyear{1992}
\cmonth{Jan}
\cyear{1992}
\cvolno{1}
\cpgs{3-28}
\title Some basic information\\
on information-based complexity theory \endtitle
\author Beresford N. Parlett\endauthor
\shortauthor{B. N. Parlett}
\shorttitle{Information-Based Complexity Theory}
\address Department of Mathematics and Computer Science 
Division of the
Electrical Engineering and Computer Science Department, 
University of
California, Berkeley, California 94720\endaddress
\date November 10, 1989 and, in revised form, July 9, 
1990\enddate
\subjclass Primary 68Q25; Secondary 65F99, 
65Y20\endsubjclass
\thanks The \kern-.39pt author \kern-.39pt gratefully
 \kern-.39pt acknowledges \kern-.39pt partial \kern-.39pt 
support
 \kern-.39pt from \kern-.39pt ONR \kern-.39pt Contract 
\kern-.39pt 
N00014-90-J-1372\endthanks
\abstract Numerical analysts might be expected to pay 
close attention to a
branch of complexity theory called information-based 
complexity theory (IBCT),
which produces an abundance of impressive results about 
the quest for
approximate solutions to mathematical problems. Why then 
do most numerical
analysts turn a cold shoulder to IBCT? Close analysis of 
two representative
papers reveals a mixture of nice new observations, error 
bounds repackaged in
new language, misdirected examples, and misleading 
theorems.
\par Some elements in the framework of IBCT, erected to 
support a rigorous yet
flexible theory, make it difficult to judge whether a 
model is off-target or
reasonably realistic. For instance, a sharp distinction 
is made between
information and algorithms restricted to this 
information. Yet the information
itself usually comes from an algorithm, so the 
distinction clouds the issues
and can lead to true but misleading inferences. Another 
troublesome aspect of
IBCT is a free parameter $F$, the class of admissible 
problem instances. By
overlooking $F$'s membership fee, the theory sometimes 
distorts the economics
of problem solving in a way reminiscent of agricultural 
subsidies.
\par The current theory's surprising results pertain only 
to unnatural
situations, and its genuinely new insights might serve us 
better if expressed
in the conventional modes of error analysis and 
approximation theory.\endabstract
\endtopmatter

\document
\ac
\heading 1. Introduction\endheading
\heading 2. High level criticisms\endheading
\heading 3. Preliminaries\endheading
\heading 4. On the optimal solution of large linear 
systems\endheading
\heading 5. Optimal solution of large eigenpair problems 
\endheading
\heading {} References\endheading
\endac

\heading1. Introduction and summary\endheading
\par In 1980 Traub and Wozniakowski published a monograph 
entitled {\it A
general theory of optimal algorithms\/}, which initiated 
a new line of
research. The subject was called analytic complexity 
theory, initially, but is
now referred to as information-based complexity theory 
(IBCT hereafter). The
August 1987 issue of Bull. Amer. Math.
Soc. (N.S.) contains a summary of recent results [Pac \& 
Wo, 1987], so
acceptance of this branch of complexity theory has been 
swift.
\par One purpose of the general theory is to provide  an 
attractive setting for
the systematic study of some of the problems that have 
engaged numerical
analysts for decades. One among the programs of IBCT is 
to determine the
minimal cost of computing an approximate solution (of 
given accuracy) over {\it
all\/} algorithms that one could use which restrict 
themselves to certain
limited information about the data. It is also of 
interest to discover any
algorithms that achieve this minimal cost or, at least, 
some cost close to it.
Pride of place in IBCT is given to the information to 
which the algorithms
are limited. By its choice of problems IBCT is 
(potentially) a branch of
complexity theory that is highly relevant to numerical 
analysts. Whenever the
minimal costs can be estimated, they provide yardsticks 
against which to
measure actual algorithms. The program seems attractive.
\par The purpose of this paper is to sound a warning 
about IBCT, but first we
point out why our task is less than straightforward and 
why our observations
have more than local interest.
\par Our task would be easier if mathematics enjoyed a 
tradition of public
debate and controversy, as exists in some other fields, 
from which some kind of
consensus could emerge  if not Truth. But too many of us 
treat every
mathematical assertion as if it must be true or false, 
and if false, then
mistaken, and if mistaken, then a symptom of some lapse 
or inaptitude. Such an
oversimplification would wrongly impugn the technical 
competency of the
founders of the IBCT and miss the point of our criticism; 
we intend to show
that parts of IBCT are true, but mistaken. Our task would 
be easier if IBCT
contained a logical flaw so egregious that its mere 
exposure rendered further
discussion unnecessary. But the flaws are subtle, subtle 
enough that they
appear only after their consequences have been analyzed 
in detail. That is why
we must focus upon specific work. We have selected two 
related papers [Tr \&
Wo, 1984; Ku, 1986] to be examined in detail in \S\S4 and 
5. They concern
matrix computations, which is the area we understand best 
of all the many areas
that IBCT covers. Besides, we believe these papers are 
typical; written in a
professional manner, their analysis is not weak, and some 
of their observations
are new and interesting in a technical way regardless of 
IBCT. We claim that
their results, most of them, are {\it seriously 
misleading\/}.
\par Since the arguments in the papers are impeccable, 
the flaw must be in the
framework. Yet the definitions are laid out plainly for 
all to see and seem to
be appropriate---a puzzling situation.
\par One source of difficulty is the redefinition of 
common terms such as
``eigenvalue problem'' (see \S2.4) or ``worst case'' (see 
\S2.5) or
``information'' (see \S4.4) or ``algorithm'' (see 
\S\S2.1, 3.2, 4.4, and 5.3).
These slight twists to conventional meanings are subtle 
enough to escape notice
but entail significant consequences. The results mislead 
because they are
remembered and discussed in terms of ordinary word usage. 
Most readers  will
not even be aware of the shifts in meaning, some of which 
are due to the
tempting but artificial distinction between information 
and algorithm.
\par Another feature of IBCT that can sometimes rob 
results of their relevance
is the presence of a free parameter: the class $F$ from 
which  worst cases at
to drawn. The cost of testing membership in $F$ is 
ignored by IBCT, so the
model loses validity whenever this cost is not negligible.
\par Some of our criticisms require very little knowledge 
of the subject
matter. These criticisms are presented in the next 
section. After that, we 
get down to details, provide some background material, 
and then examine each
paper in turn. In our summaries we try to be fair, but we 
encourage the
interested reader to compare our effort with the original 
work.
\par A handful of reservations about IBCT have appeared 
in print. In a review
of the second monograph [Tr, Wo \& Wa, 1983], Shub [Shu, 
1987] gives a couple of
instances of unnatural measures of cost. In [Ba, 1987], 
Babuska calls on
researchers in IBCT to make their models more realistic. 
We concur but note
that the model may be so flexible as to embrace pointless 
investigations as
readily as pertinent ones.
\par We make no complaint that IBCT ignores the roundoff 
error that afflicts
implementation on digital computers. First, a good 
understanding of
computation in exact arithmetic is a prerequiste for 
tackling practical issues.
Second, we must acknowledge that a large part of 
theoretical numerical
analysis confines itself to the comforts of exact 
arithmetic.
\par IBCT has already produced a large body of results, 
some of them surprising
and consequently of potential interest. Yet each 
surprising result known to us,
in worst-case analysis, holds only within a model 
sufficiently unnatural as to
forfeit attention from numerical analysts. This is a pity 
because IBCT
certainly permits realistic models, and there is plenty 
to do; the
investigation of average case complexity of approximately 
solved problems is in
its infancy. It would take only a few illuminating 
results concerning some
reasonable models to restore faith in the program 
launched by the general
theory of optimal algorithms. So, it is a pleasure to 
acknowledge the recent
and interesting work by Traub and Wozniakowski on the 
average case behavior o
the symmetric Lanczos algorithm [Tr \& Wo, 1990]. 
However, we must note that the
infrastructure of IBCT plays little role in this analysis.
\par The incentive to write this essay came from 
discussions held during the
workshop at the Mathematical Sciences Research Institute 
at Berkeley in
January 1986 under the title Problems Relating Computer 
Science to Numerical
Analysis.
\heading 2. High level criticisms\endheading
\subheading{{\rm 2.1.}\quad This is not complexity 
theory} Numerical analysis
and complexity theory are palpably different subjects. 
Complexity theory (CT
hereafter) seeks to determined the {\it intrinsic\/} 
difficulty of certain
tasks; whereas much of theoretical numerical analysis (NA 
hereafter) has been
concerned with studying classes of {\it algorithms\/}, 
analyzing convergence
and stability, developing error bounds (either a priori 
or a posteriori), and
detecting either optimal or desirable members of a class 
according to various
criteria. Clearly CT has more ambitious goals than does NA.
\par One major theme of IBCT is to find the minimal cost 
of achieving a
certain level of approximation for the hardest case in a 
given problem class
$F$, using any algorithm that confines itself to certain 
partial information
about the case. One of the papers we examine is concerned 
with the solution of
large systems of linear equations and the other with the 
matrix eigenvalue
problem (see [Tr \& Wo, 1984; Ku, 1986]). Now we can 
formulate our first
complaint, one that applies to the results in both papers.
\ext
The theorems say nothing about the intrinsic cost of 
computing an approximate
solution to either of the problems mentioned above 
because the specified
information is not naturally associated with the task but 
is acquired when a
certain class of numerical methods is employed.
\endext
\par The class is sometimes called the Krylov subspace 
methods; one is {\it
not\/} given a matrix $A$ explicitly, but instead a few 
vectors $b$, $Ab$,
$A^2b$, $A^3b,\dotsc, A^jb$, and one wishes to 
approximate the solution of a
system of linear equations $Ax=b$ or some specified 
eigenpair of $A$. More
details are given in \S\S3.2 and 3.3. So the invitation 
to minimize cost over
{\it all\/} algorithms subject to the given information 
turns out, in these
cases, to amount to the quest for the best that can be 
achieved at {\it each
step\/} of a Krylov subspace method. This is exactly the 
sort of work that
numerical analysts do.
\par We do not wish to belabor the obvious, but our 
suggestion that, in the
these cases, IBCT has the appearance of CT without the 
substance is important
for the following reason. It might be claimed that we 
interpret IBCT results as
though they were results about Krylov  subspace methods 
(i.e., NA) when, in
fact, they are CT results concerning Krylov {\it 
information\/}. In other
words, perhaps we are guilty of looking at the world 
through NA spectacles and
missing that subtle difference of emphasis characteristic 
of CT. This
possibility needs to be considered, but the stubborn fact 
remains that
restricting information to Krylov information is not part 
of the linear
equations problem nor of the eigenvalue problem.
\subheading{{\rm 2.2.}\quad Free information in the 
problem class} \kern-.39pt The
 \kern-.39pt ingredient \kern-.39pt of \kern-.39pt IBCT 
that allows it to generate irrelevant results is the 
problem
class $F$ (see the second paragraph in \S2.1). $F$ did 
not appear in our brief
description of the theory in the second paragraph of \S1 
because it is not a
logically essential ingredient but rather a parameter 
within IBCT. Let us
describe the role it plays. There is a task $T$ to be 
accomplished; there is an
information sequence $N=(N_1, N_2,N_3,\ldots)$ coupled 
with a measure of cost
$(N_j$ costs $j$ units); and there is $F$. For each $N_j$ 
the only algorithms
admitted for consideration are those that restrict 
themselves to use $N_j$ and
standard auxiliary computations. For a worst-case 
complexity analysis, the
main technical goal is to determine the minimal cost, 
over admissible
algorithms, required to achieve $T$ for the most 
difficult problem {\it
within\/} $F$ that is consistent with $N$. This minimal 
cost may be called
$C(F,N,T)$.
\par Suppose now that $F_1\!\subset \!F_2$ and that 
$C(F_1\!,N\!,T)\!<\!C(F_2\!,N\!,T)\!$. Such a
result is of little relevance to the achievement of $T$ 
unless one can
determine that a problem lies within $F_1$ rather than 
$F_2$. To put the matter
in other words, we might say that knowledge of membership 
in $F$ is information
and should have a cost attached to it. Whenever $F$ is 
very large (for example,
the class of continuous functions or the class of 
invertible matrices), it is
realistic to assign no cost to it. On the other hand, 
there are examples (see
\S4.4) where it may be as expensive to ascertain 
membership in $F$ as to
achieve $T$, given $N$, over a larger class of problems. 
In such cases
$C(F,N,T)$ bounds one part of the expense while ignoring 
the other. Let
$C(N,T)$ denote $C(F,N,T)$ when $F$ is as large as
possible.
\par We may reformulate the minimax quantity $C(N,T)$ 
with the aid of a useful
piece of notation. To each $N_j$ there is a set 
$\widehat{V}_j$ of problems
(matrices in our case) that are indistinguishable by 
$N_j$. The set $\widehat{V}
_j$, $j=1,2,\dotsc, n$, are nested and eventually reduce 
to a singleton.
Associated with any approximation $z$ is the set $R_j(z)$ 
of indistinguishable
residuals (e.g., $R_j(z)=b-\widehat{A}z$, $\widehat{A}\in 
\widehat{V}_j$, for
the linear equations problem $Ax=b)$. The goal is to find 
the smallest natural
number $k$ such that there is a $z$ for which $R_k(z)$ 
lies in the target ball
(e.g., $B(0,\varepsilon \|b\|)$, the ball in $R^n$ 
centered at the origin with
radius $\varepsilon \|b\|)$. This is $C(N,T)$.
\par This formulation reveals several things. First, the 
admissible algorithms
cited in the minimax formulation of $C(N,T)$ are not 
really needed; what
matters is the size of $R_k(z)$ for various $z$. Second, 
one reason why there
is very little in the NA literature on the problem of 
finding the minimal $k$
is that for most interesting tasks $k=n$, the sets 
$R_k(z)$ are just too big, so
the problem is not interesting.
\par One way to reduce the indistinguishable sets is to 
introduce a subclass
$F$ and to use $V_j=\widehat{V}_j\cap F$ in place of 
$\widehat{V}_j$. This was
discussed above. For approximation theory there is no 
objection to the
introduction of unknown quantities that might 
characterize $F$. However, as
mentioned above. IBCT seems to use $F$ as a tuning 
parameter designed to keep
$k<n$.
\subheading{{\rm 2.3.}\quad Spurious challenges} The 
optimality properties of
various Krylov subspace methods are well known (see [Sti, 
1958]). IBCT's claim
to have something new to add is based on the suggestion 
that its theory
considers any algorithm (confined to Krylov information) 
including those that
give approximations $z$ from {\it outside\/} the Krylov 
subspace. See the
quotation in \S4.1. The trouble with this apparent 
novelty is that it is not
possible to evaluate the residual norm $\|b-Az\|$ for 
these external $z$
because there is no known matrix $A$ (only Krylov 
information). So how can an
algorithm that produces $z$ verify whether it has 
achieved its goal of making
$\|b-Az\|<\varepsilon \|b\|$? Perhaps that is why no such 
new algorithm is
actually exhibited? IBCT's suggestion that it goes beyond 
the well-known
polynomial class of algorithms is more apparent than real.
\subheading{{\rm 2.4.}\quad A new eigenvalue problem} The 
task of computing
some or all the eigenvalues of a matrix is acknowledged 
to be of practical
importance. When only a few eigenvalues of a large order 
matrix are wanted, one
seeks either the smallest eigenvalues, or the largest, or 
all in a given
region.
Unfortunately, [Ku, 1986] makes a subtle change in the 
problem. The redefined
goal asks for {\it any\/} approximate {\it eigenpair\/} 
(value $\lambda $ and
vector $x)$ without reference to where in the spectrum 
the approximate
eigenvalue may lie. Of course a theorist is entitled to 
investigate any
problem he or she chooses. However, we have yet to hear 
of {\it any\/} use for
such output. Our complaint is that {\it no indication\/} 
is given that the goal
is an unusual one. Very few readers would realize that 
the familiar relevant
eigenvalue problem has been bypassed. Indeed, we missed 
the point ourselves
until a friend pointed it out.
\par It is standard practice to use the size of the 
residual norm
$(\|Ax-x\lambda \|)$ as the means by which to decide 
whether a {\it
specified\/} approximate eigenvalue is accurate enough. 
IBCT twists things
around and makes a small residual norm into the goal. It 
is the old
philosophical error of mistaking the means for the end.
\subheading{{\rm 2.5.}\quad A confusion of worst cases} 
An important feature of
Krylov information $\{b,Ab,A^2b,\ldots\}$ (see the third 
paragraph in \S2.1) is
the so-called starting vector $b$ which is, of course, 
quite independent of the
goal of computing eigenvalues. There are two different 
factors that can
increase the cost of using this information to 
approximate eigenvalues of $A$.
One is an intrinsic difficulty; some matrices in the 
given class may have
unfortunate eigenvalue distributions. The other is that 
$b$ may be poorly
chosen. Instead of separating the effects of these 
factors, the eigenvalue
paper combines them and ends up analyzing Krylov 
information with a worst
possible starting vector even though satisfactory staring 
vectors are easy to
obtain. The fact that $b$ is treated as prescribed data 
is quite difficult to
spot. This situation is in stark contrast to the linear 
equations problem where
$b$ is part of the problem of solving $Ax=b$.
\par The study of worst choices for $b$ is not without 
interest (see [Sc, 1979],
for example). Such studies are relevant to the inverse 
eigenvalue problem but
not to the complexity of approximating eigenvalues via 
Krylov subspaces.
\par Section 5 discusses the issue in more detail, but 
the conclusion is that
the wrong placing of $b$ twists the model away from its 
original target. It is
only this distorted model that permits the proof of the 
surprising results
in [Ku, 1986, Abstract].
\heading 3. Preliminaries\endheading
\subheading{{\rm 3.1.}\quad A word on matrix 
computations} \kern-.39pt The \kern-.39pt subject 
\kern-.39pt begins
\kern-.39pt with \kern-.39pt three basic tasks.
\roster
\item"(i)" Solve systems of linear algebraic equations; 
\kern-.39pt 
written as \kern-.39pt
 $Ax\!=b$ with
$x$, $b$ column vectors and $A$ a matrix.
\item"(ii)" Compute eigenvalues and eigenvectors of $A$.
\item"(iii)" Solve least squares problems, i.e. minimize 
$\|b-Ax\|_2$ over all
$x$.
\endroster
\par There are very satisfactory programs for 
accomplishing these tasks 
when the matrices are small. An $n\times n$ matrix is 
said to be {\it small\/}
if a couple of $n\times n$ arrays can be comfortably 
stored in the fast memory
of the computer. These days a $50\times 50$ matrix would 
be considered small on
most computer systems. The reader may consult [Pa, 1984] 
for more information.
The methods in use make explicit transformations on the 
given matrix. There are
one or two open problems concerning convergence of some 
methods, but by and
large the small matrix problem is in satisfactory 
condition with respect to
conventional one-calculation-at-a-time (sequential) 
computers.
\par One reason for introducing this slice of history 
into the discussion is to
bring out the fact that computation with large order 
matrices (say, $5000\times
5000)$ is a somewhat different game from computation with 
small ones. Sometimes
the very problem itself changes. For tasks (i) and (iii) 
the goal does remain
the same but often the product $A\nu$, for any vector 
$\nu$, can be formed cheaply,
so one seeks methods that exploit this feature and do not 
factor $A$. For the
eigenvalue problem, there are many applications where 
only a few eigenvalues
are wanted, perhaps 30 of 5000, and it is desirable to 
avoid changing $A$ at
all. Thus the task has changed; there is no desire to 
diagonalize $A$.
\par For all three problems it often happens that a 
sequence of similar cases
must be treated as parameters are changed. This leads to 
an updating problem
and ends our brief historical digressions.
\par As usual $R$ denotes the real numbers, $C$ the 
complex numbers, and $R^n$
is the vector space of real $n$-tuples.
\subheading{{\rm 3.2.}\quad A word on information-based 
complexity} We describe
a simple version of IBCT that is used in the two papers 
to be examined. It does
not use the full panoply of concepts in the monograph or 
its sequel [Tr, Wo \&
Wa, 1983].
\par There are a few essential ingredients that are best 
seen in a specific
context.
\roster
\item"(1)" A class $F$, e.g., $F=\{B\: B\in R^{n\times 
n}$ symmetric, positive
definite\};
\item"(2)" A task, e.g., given $b\not=0$ in $R^n$, and 
$\varepsilon >0$, find
$x$ in $R^n$ s.t. $\|Ax-b\|<\varepsilon \|b\|$, $A\in F$. 
Here $\|\boldcdot\|$
is some given norm.
\endroster
\par For the eigenvalue problem, the task is stated as: 
find $x\in C^n$ and
$\rho\in C$ such that 
$\|\widetilde{A}x=x\rho\|\le\varepsilon $ for all
$\widetilde{A}$ in $F$ that are indistinguishable from 
$A$. The defects in this
definition were mentioned in \S2.5.
\roster
\item"(3)" information $N=(N_0,N_1,\ldots)$, e.g., 
$N_j(A,b)=\{b,Ab,\ldots,
A^jb\}$ for natural numbers $j\le n$, $A\in F$;
\item"(4)" A measure of cost, e.g., $j$ units for $N_j$. 
In this model the
forming of linear combinations of vectors is free.
\endroster
Items (2) and (3) do not make it clear that $A$ is not 
known explicitly. There
is more discussion of this point in \S4.3.
\par To use an algorithm, restricted to the information 
$N$, in order to solve
a problem in $F$ will entail cost that may vary with the 
problem.
 The primary goal of worst-case 
IBCT is to minimize, over all such algorithms, the 
maximum, over all
problems in $F$, of that cost. Determining the
minimum, even roughly, is worthwhile even if an algorithm 
that achieves
the minimum is not known. There is an analogous 
average-case theory.
\par This certainly has an appeal.
\par Please note that, in our example, (3) puts this 
theory firmly within
numerical analysis. This is because the information in 
this example, and it is
typical, is {\it not\/} part of the task. The information 
$N_j$ will only be
available if a certain type of method is invoked. 
Consequently the theory is
{\it not\/} addressing the intrinsic cost, or difficulty, 
of solving linear
systems but is confined to seeking the number of needed 
steps within a chosen
class of methods. This is what numerical analysts do, and 
have done, from the
beginning.
\par In the 1970s the word ``complexity'' was reserved 
for the intrinsic
difficulty of a task and the word ``cost'' was used in 
connection with a
particular algorithm. For example, see [Bo \& Mu, 1975] 
and [Wi, 1980].
However, it is now common to talk about the complexity of 
an algorithm as a
synonym for cost. This extension of the term complexity 
does no great harm. What
is misleading is that the notion of information {\it 
appears\/} to be
independent of any algorithm. This allows the theory to 
talk about the set of
all algorithms that confine themselves to the given 
information. As indicated
in the previous paragraph, this way of talking may 
sometimes be a reformulation
of the standard practice of optimizing over a specified 
family of methods.
\par For $j<n$ the information $N_j(A,B)$ is partial; 
there are many matrices
that are indistinguishable from $A$ in the sense that 
each of them generates
the same set of $j+1$ vectors. The basic technical 
concept in the theory, the
{\it matrix index\/} $k(\Phi,A)$ of an algorithm $\Phi$, 
is the minimal cost
using $\Phi$ to guarantee achievement of the task for all 
matrices in $F$ that
are indistinguishable from $A$. There is more discussion 
in \S4.2 and \S2.2.
\par The theory seeks $\operatorname{min}_\Phi k(\Phi,A)$ 
and other related
quantities while $\Phi$ is restricted to information $N$. 
This minimum is the
complexity of the task.
\par IBCT is careful to distinguish itself from the older 
speciality now called
arithmetic complexity, which is concerned with discrete 
problems such as the
minimal number of basic arithmetic operations required 
 to form the product of
two $n\times n$ matrices (see, for example, [Ra, 1972; 
Str, 1969; Sch \& Str,
1971; Wi, 1970]). Another way to model some aspects of 
scientific computing was
introduced in 1989 by Blum, Shub, and Smale, [Bl, Sh \& 
Sm, 1989]. Their
algebraic complexity has the cost of a basic arithmetic 
operation independent
 of the operands
just as it does in scientific computation. They analyze 
computation
over rings other than $Z$ and so include $R$ and $C$.
\subheading{{\rm 3.3.}\quad Krylov subspaces} Here we 
sketch the conventional
wisdom on this topic. These subspaces of $R^n$ are 
defined by
$$K^j=K^j(A,b)=\operatorname{span}(b,Ab,\dotsc, 
A^{j-1}b).$$
There is no loss in assuming that 
$\operatorname{dim}K^j=j$. Information $N_j$
permits computation of any vector $v$ in $K^{j+1}$, not 
just $K^j$, at no
further cost provided that the coefficient $\gamma_i$ in
$v=\sum^j_{i=0}\gamma_i(A^ib)$ are known.
\par On the practical side the dominant question for the 
experts has been how
to obtain a computationally satisfactory basis for $K^j$. 
Round off error
destroys the expected linear independence of the computed 
vectors. Some
researchers maintain that it is more efficient to use a 
highly redundant
spanning set rather than a basis. Others recommend the 
additional expense of
computing an orthonormal basis. In any case, it is the 
computation of the basis
or spanning set, along with multiplication of vectors by 
$A$, that is the
expensive part of the computation. The model of cost used 
in IBCT theory
reflects this quite well. It is the number of {\it 
steps\/} that matters. We
think of a step as adding one more vector to the Krylov 
sequence
$\{b,Ab,A^2b,\ldots\}$.
\par One of the founders of Krylov space methods, C. 
Lanczos [La, 1952],
proposed a useful basis for $K^j$ with the property that 
the projection of
symmetric $A$ onto $K^j$ is a symmetric tridiagonal 
$j\times j$ matrix $T_j$.
Tridiagonal matrices are easy to handle. With a basis in 
hand, there are diverse
tasks that can be accomplished. Here is a partial list.
\roster
\item"(i)" Compute an approximation $x^{(j)}$ to 
$A^{-1}b$ such that
$x^{(j)}\in K^j$ and its residual $b-Ax^{(j)}$ is 
orthogonal to $K^j$. It so
happens that $\|b-Ax^{(j)}\|$ may be computed without 
forming $x^{(j)}$ so
there is no need to compute unsatisfactory $x^{(j)}$. When
$A\in\operatorname{SPD}$ (symmetric positive definite 
matrices) then $x^{(j)}$
coincides with the output of the conjugate gradient 
algorithm.
\item"(ii)" Compute the vector $u^{(j)}$ that minimizes 
$\|b-Av\|$ over all
$v\in K^j$ (not $K^{j+1})$. This is the MR (minimum 
residual) algorithm. The
extra vector $A^jb$ is needed to ascertain the 
coefficients in the expansion of
$u^{(j)}$.
\item"(iii)" Compute some, or all, of the Rayleigh-Ritz 
approximations
$(\theta_i,y_i)$, $i=1,\dotsc, j$, to eigenpairs of 
symmetric $A$. Here
$\theta_i\in R$ and $\{y_1,\dotsc, y_j\}$ is an 
orthonormal basis for $K^j$. For
each $i$, $Ay_i-y_i\theta_i$ is orthogonal to $K^j$.
\endroster
\par Krylov subspace methods are not really iterative. 
All the basic tasks
mentioned in \S3.1 are solved exactly (in exact 
arithmetic) in at most $n$
steps. However the interest in this approach is due to 
the fact that in many
instances far fewer than $n$ steps are required to 
produce acceptable
approximations. In other words, to take $n$ steps is the 
practical
equivalent of failure. However, for each basic task there 
are 
data pairs $(A,b)$ that do require $n$ steps even for 
loose tolerances such as
$\varepsilon =10^{-3}$, so research has focussed on 
explanations of why so
often the cost is much smaller. The gradual realization 
of the efficacy of
changing the $Ax=b$ problem to an equivalent one via a 
technique called 
preconditioning has enhanced the use of Krylov subspace 
methods.
\par In a sense the {\it convergence\/} of all these 
methods is completely
understood in the symmetric case and is embodied in the 
error bounds published
by Kaniel and improved by Paige and Saad (see \cite{Ka, 
1966; Sa, 1980; Pa,
1980} for the details). The error depends on two things; 
the eigenvalue
distribution of $A$ and the components of the starting 
vector $b$ along the
eigenvectors. Of course, all this analysis supposes a 
given matrix $A$, not a
set of indistinguishable matrices.
\par From this conventional viewpoint the thrust of these 
two complexity papers
is to see to what extent the standard algorithms (CG, MR, 
Lanczos) do not make
best use of the information on hand. Recall that $N_j$ 
contains an extra vector
not in $K^j$. This is a reasonable project and the 
results can be expressed
within the usual framework. The term ``complexity 
theory'' appears to a
numerical analyst like window dressing.
\heading 4. On the optimal solution of large linear 
systems\endheading
\par These sections offer a description of and commentary 
on [Tr \& Wo, 1984].
\subheading{{\rm 4.1.}\quad Spurious generality} Here is 
a quotation from the
introduction:
\ext
We contrast our approach with that which is typical in 
the approximate solution
of large linear systems. One constructs an algorithm 
$\Phi$ that generates a
sequence $\{x_k\}$ approximating the solution $\alpha 
=A^{-1}b$; the
calculation of $x_k$ requires $k$ matrix-vector 
multiplication and $x_k$ lies
in the Krylov subspace spanned by $b\!, Ab\!,\dotsc\!, 
A^kb$. The algorithm $\Phi$
is often chosen to guarantee good approximation 
properties of the sequence
$\{x_k\}$. In some cases $\Phi$ is defined to minimize 
some measure of the
error in a {\it restrictive} class of algorithms. For 
instance, let this class be
defined as the class of `polynomial' algorithms; that is
$$\alpha -x_k=W_k(A)\alpha ,\quad\text{where }W_k(0)=1.$$
Here $W_k$ is a polynomial of degree at most $k$.
$$\ldots$$
[Some \kern-.39pt omitted \kern-.39pt sentences 
\kern-.39pt define \kern-.39pt the \kern-.39pt minimum 
\kern-.39pt 
re\-sidual and conjugate gradient
algorithms.]
\newline
It seems to us that this procedure is unnecessarily 
restrictive. It is not
clear, a priori, why an algorithm has to construct $x_k$ 
of the form $\alpha
-x_k=W_k(A)\alpha $. Indeed, we show that for 
orthogonally invariant classes of
matrices the minimum residual algorithm (MR) is within at 
most one matrix
vector multiplication of the lower bound without any 
restriction on the class
of algorithms. However, if the class is not orthogonally 
invariant, the
optimality of MR may disappear.
\endext
\par Our first point was made earlier. The information 
$N$ does not come with
the linear equations problem. The brief answer to the 
quoted rhetorical
question (why must an algorithm construct $x_k$ of the 
given form?) that serves
to justify the whole paper is the following. To any 
vector $x$ {\it not\/} in
the Krylov subspace $K^k$, there is an admissible matrix 
$A$, such that the
residual norm is as large as you please. This holds even 
when $A$ is required
to be symmetric and positive definite. An admissible 
matrix $A$ is one that is
consistent with the Krylov information (more on this 
below).
\subheading{{\rm 4.2.}\quad Definitions and optimality} 
In this section we put
down the definitions made in [Tr \& Wo, 1984]. Our 
comments are reserved for
the next section.
\roster
\item"(i)" Let $F$ be a subclass of the class 
$\operatorname{GL}(n,R)$ of
$n\times n$ nonsingular real matrices.
\item"(ii)" Let $b\in R^n$ with $\|b\|=(b,b)^{1/2}=1$ be 
given. For
$0\le\varepsilon <1$, find $x\in R^n$ such that 
$\|b-Ax\|\le\varepsilon $ (it
would have been kinder to add $A\in F)$.
\item"(iii)" Krylov information: $N_j(A,b)=\{b,Ab,\dotsc, 
A^jb\}$,
$j=0,1,\ldots$.
\item"(iv)" Measure of cost: $N_j$ costs $j$ units.
\item"(v)" An algorithm $\Phi=\{\phi_j\}$ is a sequence 
of mappings $\phi_j\:
N_j(F,R^n)\rightarrow R^n$.
\item"(vi)" The set of indistinguishable matrices for 
given $N_j(A,b)$:
$$V(N_j(A,b))=\{\widetilde{A}:\widetilde{A}\in F\:
N_j(\widetilde{A},b)=N_j(A,b)\}.$$
\item"(vii)" The matrix index of an algorithm $\Phi$:
$$\hphantom{k}k(\Phi,A)=\operatorname{min}\left\{j: 
\max_{\tilde{A}\in
V(N_j)}\|b-\widetilde{A}x_j\|\le\varepsilon 
\right\},\qquad N_j=N_j(A,b),$$
where $x_j=\Phi_j(N_j(A,b))$. If the set of $j$ values is 
empty, then
$k(\Phi,A)=\infty$.
\item"(viii)" The class index of an algorithm $$\Phi\: 
k(\Phi,F)=
\operatorname{max}_{B\in F}k(\Phi,B).$$
\item"(ix)" The optimal matrix index: 
$k(A)=\operatorname{min}_\Phi k(\Phi,A)$
over $\Phi$ restricted to $N$. 
\item"(x)" The optimal class index: 
$k(F)=\operatorname{max}_{B\in F}
k(B)$.
\item"(xi)" Strong optimality: $\Phi$ is strongly optimal 
iff $k(\Phi,B)=k(B)$,
for each $B\in F$.
\item"(xii)" Optimality: $\Phi$ is optimal iff 
$k(\Phi,F)=k(F)$.
\item"(xiii)" Almost strong optimality: $\Phi$ is almost 
strongly optimal iff
$k(\Phi,B)\le k(B)+c$, for every $B\in F$, for some small 
integer $c$.
\endroster
\rem{Remark \RM1} Since $A^ib=A(A^{i-1}b)$, it follows 
that Krylov information
$N_j(A,b)$ requires $j$ applications of the operator $A$. 
That is why the cost
is given as $j$ units. In practice, one uses better bases 
for the Krylov
subspace $K^j$ than is provided by $N_j(A,b)$; but for 
theoretical purposes,
this modification may be ignored.
\endrem
\rem{Remark \RM2} It can happen that $k(A)\ll k(F)$. For 
this reason, it is of
interest to find algorithms with small matrix index.
\endrem
\rem{Remark \RM3} For simplicity, the dependence of all 
concepts on $n$, $N_j$,
$b$, and $\varepsilon $ is suppressed. The idea is to 
compute $k(A)$ and $k(F)$
for interesting classes $F$ and to find strongly optimal 
or optimal algorithms
if possible.
\endrem
\subheading{{\rm 4.3.}\quad Discussion of the basic 
concepts} In \S2 we pointed
out how misleading it can be to compute complexity for 
restricted classes $F$,
which are difficult to discern in practice. Here we wish 
to point out that $F$
is introduced into the basic definitions, such as $V$ in 
(vi), and there is no
need for it.
\par To add to any confusion, the basic definitions do 
not make clear the role
of $A$. In the context of numerical analysis there is a 
particular matrix $A$
on hand and this permits one to test the residual 
$r=b-Av$ for any vector $v$.
However, in the context of IBCT, that is not quite the 
case. In this game we
consider a specific $A$, but it is not available 
explicitly. That odd
situation certainly warrants some discussion, and it 
faithfully reflects the
state of affairs at a typical step in a Krylov subspace 
method. The matrix is
hidden inside a subprogram, and the user has only a basis 
for the Krylov
subspace corresponding to Krylov information 
$N_j(A,b)=\{b,Ab,\dotsc, A^jb\}$.
Associated with $N_j(A,b)$ is
$$\widehat{V}(N_j(A,b))=\{\widetilde{A}:N_j(%
\widetilde{A},b)=N_j(A,b)\},$$
the set of matrices indistinguishable from $A$ by 
$N_j(A,b)$. Contrast
$\widehat{V}$ with $V$ in \S4.2 (vi).
\par With $\widehat{V}$ defined above, the natural 
definition of the matrix
index of an algorithm $\Phi$ is 
$$\hat k(
\Phi,A)=\operatorname{min}\left\{j:\max_{\widetilde{A}\in%
\widehat{V}(N_j)}\|b-
\widetilde{A}x_j\|\le\varepsilon \right\},\qquad 
N_j=N_j(A,b),$$
where
$$x_j=\Phi_j(N_j(A,b)).$$
If the set of $j$ values is empty then $\hat 
k(\Phi,A)=\infty$. Please note
that in contrast to (vii) in \S4.2 there are no hidden 
parameters in $\hat k$.
It is the first step in the process at which the task is 
accomplished by
$\Phi$ for all matrices indistinguishable from $A$ by 
$N_1$ through $N_{\hat
k}$. Then the optimal index is
$$\hat k(A)=\min_\Phi\hat k(\Phi,A)$$
over all $\Phi$ such that $\Phi_j$ uses only $N_j(A,b)$ 
and standard arithmetic
operations.
\par There is no logical need for $F$. However, given a 
class $F$, one may
define
$$\hat k(\Phi,F)=\max_{B\in F}\hat k(\Phi,B);\qquad \hat 
k(F)=\min
_\Phi\hat k(\Phi,F).$$
\par Why did IBCT not follow this simple approach? Why 
does IBCT use $V=
\widehat{V}\cap F$ to define the matrix index $k(\Phi,A)$ 
and thus suppress the
role of $F$? The reason, we suspect, is that with these 
natural definitions the
``polynomial'' algorithms, deemed to be restrictive in 
the introduction to [Tr
\& Wo, 1984] are mandatory; and, consequently, IBCT has 
nothing new to offer.
Here is a result of ours that shows why the nonpolynomial 
algorithms are of no
interest in worst case complexity, i.e., $\hat 
k(\Phi,A)=\infty$.
\proclaim{Theorem} \kern-.39pt Assume \kern-.39pt 
$A\!=\!A^t\!\in\! R^{n\times n}\!$. \kern-.39pt Let
$K^j\!=\!\operatorname{span}(b\!,Ab\!,\dotsc\!, 
A^{j-1}\!b)$ have dimension $j$ $(<n)$ and
$\|b\|=1$. To each $v\notin K^j$ there exists 
$\widetilde{A}\in\widehat{V}
(N_j(a,b))$ such that $\|b-\widetilde{A}v\|>1$.
\endproclaim
\demo{Sketch of proof} (1) Should $b$ happen to be 
orthogonal to some
eigenvectors \kern-.39pt of \kern-.39pt $A$, \kern-.39pt 
it \kern-.39pt is \kern-.39pt always \kern-.39pt possible 
\kern-.39pt to \kern-.39pt choose \kern-.39pt an 
\kern-.39pt $\overline{A}
\!\in\!\widehat{V}(N_j(A\!,b))$ such that $b$ is not 
orthogonal to any of $
\overline{A}$'s eigenvectors. If necessary replace $A$ by 
$\overline{A}$.
\roster
\item"(2)" There is a distinguished orthonormal basis for 
$K^j$ that can be
extended to a basis for $R^n$ and in which $A$ is 
represented by the matrix
$$\bmatrix
T&E^t\\
E&U\endbmatrix$$
where $T=T^t\in R^{j\times j}$, $E$ is null except for a 
$\beta \not=0$ in the
top right entry, $U$ is unknown. Moreover $T$ and $\beta 
$ are determined by 
$N_j(A,b)$.
\item"(3)" In this distinguished basis $b$ is represented 
by
$e_1=(1,0,0,\dotsc, 0)^t$ and let any $v\in K^j$ be 
represented in partitioned
form $(f,g)$ where $f\in R^j$, $g\in R^{n-j}$. By 
hypothesis $g\not=0$, since
$v\notin K^j$, and
$$\|b-\widetilde{A}v\|^2=\|e_1-Tf-E^tg\|^2+\|Ef+
\widetilde{U}g\|^2,$$
where $\widetilde{U}$ is the $(2,2)$ block in the 
representation of
$\widetilde{A}$ and is uncontrained.
\item"(4)" To each $\widetilde{U}\in 
R^{(n-j)\times(n-j)}$ there is an
$\widetilde{A}\in \widehat{V}(N_j(A,b))$ and for any 
$g\not=0$, there exists
$\widetilde{U}$ such that

\vskip 6pt
$$\|Ef+\widetilde{U}g\|^2=\|e_1\beta f(j)+
\widetilde{U}g\|^2>1.$$
\endroster
\par\vskip 6pt In particular, it is possible to select 
$\widetilde{U}$ to be symmetric,
positive definite, and diagonal. Here ends the sketch of 
the proof.
\enddemo
\par Symmetry is not needed in the result given above, If 
$A$ is not symmetric
there is still a distinguished orthonormal basis for 
$K^j(A,b)$ and $R^n$ such
that $b$ is represented by $e_1$, and $A$ is represented by

\vskip 6pt
$$\bmatrix H&J\\
E&L\endbmatrix.
$$

\vskip 6pt
\noindent Now $H$, $E$, and the first
column of $J$ are determined by $N_j(A,b)$. Moreover,
$Je_1\not=0$ and for $\widetilde{A}$ indistinguishable 
from $A$ we have
$$\|b-\widetilde{A}v\|^2=\|e_1-Hf-Jg\|^2+\|e_1\beta f(j)+
\widetilde{L}g\|^2.$$
This can be made as large as desired. In the language of 
IBCT, $\hat k
(\Phi,A)=\infty$ for any $\Phi$ such that 
$\Phi(N_j(A,b))$ takes values outside
$K^j(A,b)$.
\par Only two choices were left to IBCT.
 Either turn away to unsolved problems or cut
down the number of indistinguishable matrices by using

\vskip 6pt
$$V=\widehat{V}\cap F$$

\vskip 6pt
\noindent instead of $\widehat{V}$.
\par Here is the quandary for IBCT. If $F$ is chosen too 
small, the model loses
realism. If $F$ is allowed to be large, then the standard 
``polynomial'' regime
is optimal.
\subheading{{\rm 4.4.}\quad Discussion of results} The 
main result of the paper
concerns the minimal residual algorithm MR: this 
polynomial algorithm's output
for the information $N_j=\{b,Ab,\dotsc, A^jb\}$ is the 
vector in the Krylov 
subspace $K^j=\operatorname{span}(b,Ab,\dotsc, A^{j-1}b)$ 
that minimizes 
the residual norm, so MR is optimal in $K^j$. Given $N_j$ 
MR needs
$k(\operatorname{MR},A)$ steps to achieve an $\varepsilon 
$-approximation in the
worst case. However, IBCT says
\proclaim{Theorem 3.1 \rm [Tr \& Wo, 1984]} If $F$ is 
orthogonally invariant
then
$$k(\operatorname{MR},A)\ge k(A)\ge 
k(\operatorname{MR},A)-1\quad\text{for any
}A\in F.$$
Furthermore, both the upper and lower bounds can be 
achieved.
\endproclaim
\par\noindent
Recall that $k(A)$ is the minimal number of steps over 
all admissible
algorithms.
\par The fact that MR is {\it not\/} always strongly 
optimal {\it for the
given information\/} appears to give substance to the 
theory. It will astonish
the numerical analyst, so let us look at the example that 
purports to show
that MR is not always optimal for the given information 
(Example 3.2 in the
paper). This class is
$$\widetilde{F}_\rho=\{A: A=I-B,\ B=B^t,\ \|B\|\le 
\rho<1\}.$$
When $\varepsilon $, $\rho$, and $n$ are specially 
related so that
$$q(\varepsilon )=\left\lfloor \frac{\operatorname{ln}((1+
(1-\varepsilon
^2)^{1/2})/\varepsilon )}{
\operatorname{ln}((1+
(1-\rho^2)^{1/2})/\rho)}\right\rfloor<n,$$
then MR is just beaten by another polynomial algorithm, 
called the Chebyshev
algorithm, because
$$k(\operatorname{Cheb},A)=q(\varepsilon 
),k(\operatorname{MR},A)=q(\varepsilon
)+1.$$
\par A word of explanation is in order. Recall that 
$A^jb\in N_j(A,b)$, but
$A^jb\notin K^j$. The MR algorithm needs $A^jb$ to 
compute the coefficients $
\gamma_i$ in
$$\operatorname{MR}(N_j(A,b))=\sum^j_{i=0}\gamma_i(A^ib).$$
This always beats the Cheb output from $K^j$. However, 
Cheb can use the
well-known three-term recurrence, based on $\rho$, to 
obtain its
equioscillation approximation from $K^{j+1}$,
not just $K^j$. With the right relation between 
$\varepsilon $, $\rho$, and
$n$, one has 
$$\|b-A\operatorname{MR}(N_{q(\varepsilon 
)})\|>\|b-A\operatorname{Cheb}
(N_{q(\varepsilon 
)})\|>\|b-A\operatorname{MR}(N_{q(\varepsilon )+1})\|.$$
\par Is it fair to compare them? The theory claims to 
compare algorithms
restricted solely to information $N_j$. So how could the 
Cheb algorithm obtain
the crucial parameter $\rho$? The answer is that $\rho$ 
is found in the
definition of the problem class $\widetilde{F}_\rho$! In 
other words, knowledge
that Cheb can use is passed along through the problem 
class, not the
information.
\par The important point we wish to make is not just that 
comparisons may not
be fair but that the results of IBCT tell us as much 
about the idiosyncrasies
of its framework as they do about the difficulty of 
various approximation
problems. With a realistic class such as SPD (sym. pos. 
def.), MR is optimal
(strongly) as it was designed to be and as is well known.
\par In more recent work \cite{Wo, 1985; Tr \& Wo, 1988}, 
the flaw mentioned
above appears to be corrected and the parameter $\rho$ is 
put into the
information explicitly. Again Cheb wins by 1 because it 
uses $\rho$ while MR
does not. However, this new clarity comes at the expense 
of realism; the Krylov
information is scrupulously priced while $\rho$ comes 
free. Yet membership in
$\widetilde{F}_\rho$ may be more difficult to ascertain 
than the approximate
solution.
\par Although the IBCT paper does not mention the 
possibility, Krylov
information may be used to obtain lower bounds on $\rho$ 
that get increasingly
accurate as the dimension of the Krylov subspace 
increases. Algorithms that
exploit knowledge of the spectrum will have good average 
behavior but there is
little room for improvement in the worst case.
\par The simple facts are well known: Chebyshev 
techniques are powerful and
users are willing to do considerable preliminary work to 
estimate parameters
such as $\rho$. It is not clear, and depends strongly on 
the data, when it is
preferable to use a weaker method such as MR that does 
not need extra
parameters. The result that MR is only almost strongly 
optimal is a striking
example of obfuscation. The framework of IBCT permits 
unnatural comparison of
algorithms.
\ch Embracing the conjugate gradient algorithm\endch In 
\S4 of their
paper the authors generalize the framework to cover other 
known methods such as 
the conjugate gradient algorithm. All that is necessary 
(for IBCT) is to
introduce a parameter $p$ into the basic task. Now an 
$\varepsilon
$-approximation to $A^{-1}b$ is redefined as any $x\in 
R^n$ that satisfies
$$\|A^p(x-A^{-1}b)\|<\varepsilon \|A^{p-1}b\|.$$
The cases $p=0$, $1/2$, $1$ are the most important, and 
when $p$ is not an
integer, it is appropriate to restrict attention to the 
symmetric, positive
definite subset (SPD) of $R^{n\times n}$. When $p=1$ we 
recover MR and the
results of \S3. To generalize MR to $p=0$, it is 
necessary to use the normal
equations of a given system. The new feature, slipped in 
without mention, is
that with $p<1$ the right-hand side of the new 
definitions is {\it not
directly computable\/}. So how does an algorithm know 
when to terminate?
Please note that approximation theory can present results 
that are not
computable without a blush because it merely exhibits 
relationships. In
computation, however, there is no virtue in attaining the 
desired accuracy if
this fact
 cannot be detected. Since IBCT defines an algorithm as a 
sequence of
mappings, it dodges the difficult task of knowing when to 
stop.
\par The fact that the Conjugate Gradient algorithm 
$(p=1/2)$ minimizes
$\|b-Ax\|_A$ at each step is well known (see \cite{Da, 
1967}). Nevertheless, in
practice, the algorithm is usually terminated by the less 
desirable condition
$\|b-Ax\|<\varepsilon \|b\|$ because the desirable $A$ 
norm is not available
(see [Go \& Va, 1984]).
\par For the reasons given in the two previous 
paragraphs, Theorem 4.2 in
[Tr \& Wo, 1984], which states that under certain 
technical conditions the
generalized MR algorithm $(0\le p\le 1)$ is almost 
strongly optimal, is a true
theorem about improper algorithms.
\subheading{{\rm 4.5.}\quad An interesting result} Recall 
that 
$$\align
N_j(A,b)=&\{b,Ab,\dotsc, A^jb],\\
K^{j+1}=&\operatorname{span}\{b,Ab,\dotsc, A^jb\},\\
\widehat{V}(N_j(A,b))=&\{\widetilde{A}:N_j(%
\widetilde{A},b)=N_j(A,b)\}.
\endalign
$$
\proclaim{Theorem \rm(reformualted by us from Theorem 3.1 
in [Tr \& Wo, 1984])}
If $y\in R^n$ yields an $\varepsilon $ residual norm,
$\|b-\widetilde{A}y\|<\varepsilon \|b\|$, for all
$\widetilde{A}\in\widehat{V}(N_j(A,b))$ then so does its 
orthogonal projection
$z$ onto $K^{j+1}$.
\endproclaim
\par We offer a simplified version of the argument in [Tr 
\& Wo, 1984].
\demo{Proof} There are five steps.
\varroster
\item"(i)" Either $y\in K^{j+1}$, and there is nothing 
more to prove, or
$y=z+w$, with $z\in K^{j+1}$, and $0\not=w$ orthogonal to 
$K^{j+1}$.
\item"(ii)" For the vector $w$ defined in (i) there is a 
unique symmetric
orthogonal matrix $H=H(w)$, called the {\it reflector 
that reverses\/} $w$. In
particular: (a) $Hx=x$, for $x$ orthogonal to $w$, (b) 
$Hw=-w$. 
\endroster
\par Define an auxiliary matrix $\widehat{A}$ by 
\roster\item"(iii)"\<$\widehat{A}=HAH\in V(N_j(A,b))$ 
since, by use of (ii)(a),
$\widehat{A}\,^ib=HA^iHb=HA^ib=A^ib$, $i=1,\dotsc, j$.
\endroster
\par Note that
$$\align
\widehat{A}y-b=&HAHy-b\\
=&HA(z-w-b,\quad\text{using (i) and (ii)(b).}\\
=&H(Az-b-Aw),\quad\text{using }Hb=b.\endalign
$$
This shows the crucial relationship
\varroster
\item"(iv)"\<$\|\widehat{A}y-b\|=\|Az-b-Aw\|$, since $H$ 
preserves norms.
\endroster
Hence
\varroster
\item"(v)"
$$\aligned
\|Az-b\|\le&\tfrac12\{\|Az-b-Aw\|+\|Az-b+
Aw\|\},\qquad\qquad\quad\ \\
\shoveright{\text{by the triangle
inequality},}\\
=&\tfrac12\{\|\widehat{A}y-b\|+\|Ay-b\|\},\quad\text{by 
(iv), and (i)},\\
\le&\varepsilon \|b\|,\quad\text{using the 
hypothesis}.\endaligned
$$
\endroster
Recall that $z$ is $y$'s projection onto $K^{j+1}$. Hence 
$\widetilde{A}
z=Az$ for all $\widetilde{A}\in\widehat{V}(N_j(A,b))$ and 
so
$$\|\widetilde{A}z-b\|=\|Az-b\|\le\varepsilon 
\|b\|,\quad\text{by
(v).\quad Q.E.D.}$$
\enddemo
\par This theorem explains why MR cannot lag more than 
one step behind any
algorithm that produces an $\varepsilon $ residual norm 
for {\it all\/}
matrices indistinguishable from $A$. For, {\it by 
definition\/}, MR
produces from $N_{j+1}(A,b)$ (note the increased 
subscript) the unique vector 
in $K^{j+1}$ that gives the smallest residual and is at 
least as good as the
vectors $y$ and $z$ defined in the proof above. But $y$ 
could be the output of
a rival algorithm.
\par Our formulation of the lemma omits any mention of 
the class $F$. Now it is
clear why the hypothesis that $F$ should be orthogonally 
invariant appears in
most of the theorems. Recall that $\widehat{V}(N_j(A,b))$ 
is {\it too big\/}.
To cut down the number of {\it indistinguishable\/} 
matrices, the theory uses
$\widehat{V}\cap F=V$. To make use of the theorem, it is 
necessary to have
$HAH\in F$ and this will be true provided that $F$ is 
orthogonally invariant.
\subheading{{\rm 4.5.}\quad Summary} \kern-.39pt One 
\kern-.39pt hidden \kern-.39pt defect \kern-.39pt in 
\kern-.39pt the \kern-.39pt framework \kern-.39pt for
\kern-.39pt discussing the MR algorithms is the far 
reaching feature that allows the family
$F$ to convey what most people would call free 
information behind the back of
the information operator $N$.
\par More disturbing than the previous defect is that we 
cannot see how any
algorithm other than the well-studied polynomial 
algorithms could know when it
had achieved an $\varepsilon $-approximation if it is 
restricted to the given
information. This gives rise to a feeling that [Tr \& Wo, 
1984] managed to
create an artificial problem where no real puzzle exists. 
The mentioned
theorems (3.1 and 4.2) reflect only the propensity of 
their general theory of
optimal algorithms for creating such situations.
\heading 5. Optimal solution of large eigenpair 
problems\endheading
\par This section offers a description of and commentary 
on [Ku, 1986]. The
paper demonstrates cleverness and clean exposition but, 
nevertheless, suffers
from design flaws; it equates different versions of a 
given algorithm and it
redefines a standard task. From the abstract:
\ext
The problem of approximation of an eigenpair of a large 
$n\times n$ matrix $A$
is considered. We study algorithms which approximate an 
eigenpair of $A$ using
the partial information on $A$ given by $b$, $Ab$, 
$A^2b,\dotsc, A^jb$, $j\le
n$, i.e., by Krylov subspaces. A new algorithm called the 
generalized minimal
residual algorithm (GMR) is analyzed. Its optimality for 
some classes of
matrices is proved. We compare the GMR algorithm with the 
widely used Lanczos
algorithm for symmetric matrices. The GMR and Lanczos 
algorithms cost
essentially the same per step and they have the same 
stability characteristics.
Since the GMR algorithm never requires more steps than 
the Lanczos algorithm,
and sometimes uses substantially fewer steps, the GMR 
algorithm seems
preferable.\newline
$\ldots$ The Fortran subroutine is also available via 
$\ldots$
\endext
\par This last phrase shows that the subject matter is 
firmly within the field
of numerical analysis. Implementation issues concerning 
GMR are described in
[Ku, 1985].
\subheading{{\rm 5.1.}\quad A subtle change of goal} Here 
are the first five
lines of the paper. ``Suppose 
we wish to find an approximation to an eigenpair of
a very large matrix $A$. That is, we wish to compute 
$(x,\rho)$, where $x$ is
an $n\times 1$ normalized vector, $\|x\|=1$, and $\rho$ 
is a complex number s.t.
$$\|Ax-x\rho\|<\varepsilon \tag1.1$$
for a given positive $\varepsilon $. Here $\|\boldcdot\|$ 
denotes the 2-norm.''
\par It is all too easy to assent to this statement of 
the problem and pass on
to the rest of the article. However it is {\it not\/} the 
normal eigenvalue
problem. We are not aware of any demand at all for the 
accomplishment of this
particular task. The users of eigenvalue programs 
(engineers, theoretical
chemists, theoretical physicists) want eigenvalues in 
{\it specified parts of
the spectrum\/}; occasionally, they want the whole 
spectrum. The main concern
of this article is with symmetric matrices; and because 
their eigenvalues are
real, the usual demands are for the leftmost $p$ 
eigenvalues (for some $p\le
n)$ or the rightmost $p$ eigenvalues or for all 
eigenvalues in a given
interval. Eigenvectors may or may not be wanted. There is 
nothing inherently
wrong with restricting attention to the rather special 
case $p=1$ and a few
articles (not cited by Kuczynski) have been devoted to it 
(see [O'L, Ste \& Va,
1979; Pa, Si \& Str, 1982] for the details).
\par To support our description of user's demands we 
refer to three
publications from different fields, [Cu \& Wi, 1985, 
Introduction; Je, 1977,
Ch. 7; Sh, 1977, \S6].
\par The consequences of leaving out a vital aspect of 
the usual task are most
serious, precisely when one seeks optimal performance.
\par One reason why it is so easy to overlook the 
omission in the problem
statement is that, for symmetric matrices, almost 
everyone does use the
residual norm $\|Ax-x\rho\|$ to judge the accuracy of an 
approximate eigenpair
$(x,\rho)$. However, it is not very interesting to 
minimize the residual norm
if that might yield a $\rho$ in the unwanted part of the 
spectrum. Now a pure 
mathematician is free to define his goal at will. What is 
regrettable is that
no hint is given to the reader that the goal is not 
standard.
\par We say more about the $\varepsilon $ appearing in 
(1.1) in \S5.5.
\par We mention one other fact which may be news to 
readers who are not much
concerned with eigenvalue problems. It suggest why the 
direction taken by
Kuczynski has not appeared in the literature before. If 
we are given a
symmetric matrix $A$ and seek a single eigenvalue (with 
or without its
eigenspace) then the wanted eigenvalue is almost certain 
to be either the
leftmost or the rightmost. Recall that the Rayleigh 
quotient of a column
vector $x\not=0$ in $R^n$ is the number $x^tAx/x^tx$. The 
extreme eigenvalues
of $A$ are the extreme values of the Rayleigh quotient 
over all possible
vectors $x$ in $R^n$. It happens that, for the given 
Krylov information $N_j$,
the Lanczos algorithm is optimal for this task in the 
strong sense that it
yields the leftmost and rightmost values of the Rayleigh 
quotient over all
vectors in the {\it available\/} space $K^j$. The last 
vector $A^jb$ in $N_j$
is needed to ascertain the extreme values over $K^j$. 
Thus the problem is
settled. It is a pity that this well-known fact was not 
mentioned.
\subheading{{\rm 5.2.}\quad Choosing a bad starting 
vector} The particular
aspect of Information-Based Complexity Theory adopted in 
the paper under review
is called worst-case complexity. It seeks bounds on the 
cost of computing
$\varepsilon $-approximations over all matrices in 
certain classes $F$ and over
{\it all\/} starting directions $b$. Theorems 3.1, 3.2, 
4.1, 5.1 (there is no
theorem 1.1 or 2.1) in [Ku, 1986] are examples. In 
particular, the theory must
cover what can happen with the {\it worst possible\/} 
starting vector. Theorem
3.1 is quoted in full in \S5.4.
\par There is nothing wrong with studying the worst case. 
Indeed it
 has already been done. [Sc, 1979] is a paper with the 
clear title {\it 
How to make the Lanczos algorithm converge slowly\/} in 
which the author gives
formulae for a starting vector that prevents any Rayleigh 
Ritz approximation
from converging until the final step! Scott's paper 
appeared in Mathematics of
Computation, the American Mathematical Society's 
principal outlet for numerical
analysis, but it is not referenced in [Ku, 1986]. The 
fact that {\it some\/},
Krylov subspaces can be very badly aligned with $A$'s 
eigenvectors does prevent
worst-case analysis from shedding much light on how 
Krylov subspaces approach
certain eigenvectors in the usual case of a random 
starting vector. That
study, of course, comes under average-case analysis and 
is ripe for attention.
\par Please note that this comment is quite independent 
of comparisons of GMR
and Lanczos. The point is this: the starting vector $b$ 
is a {\it free
parameter\/} in the eigenvalue problem (in contrast to 
the linear equations
problem $Ax=b)$. It is not given and may be chosen to 
{\it improve\/}
performance. In the absence of extra information, it is 
the almost universal
habit to pick $b$ with the aid of a random number 
generator. Recent theoretical
work on Lanczos has been concerned to explain why this 
choice is so powerful
(see [Sa, 1980; Pa, 1980]). Note that two quite different 
situations have been
pushed together under the label `worst case'. It is quite 
normal to consider
the most difficult matrices $A$ because they are part of 
the problem. On the
other hand, a bad $b$ is a self-inflicted handicap rather 
than a genuine
difficulty. It is the confounding of these cases that is 
unfortunate, not their 
study.
\par Returning to the eigenvalue problem, we can rephrase 
our complaint as
follows: Kuczynski's focus, perhaps unwittingly, is on {\it
Krylov-subspaces-with-worst-possible-starting-vectors}. 
What a pity that this
was not emphasized! The numerical  examples given in the 
paper are not without
interest. The starting vector there, though not perhaps 
worst possible, is very
bad. Both methods, GMR and Lanczos converge very slowly. 
The chosen matrices
are extensions of the one used by Scott to illustrate how 
bad Rayleigh Ritz
approximations can be (see \cite{Pa, 1980, p. 218}).
\par We ran these examples with our local Lanczos 
program. It uses a random
starting vector, and convergence was quite satisfactory. 
The results are given
in \S5.5.
\par There is a different context in which the focus on 
worst starting values is
much more relevant. The GMR algorithm presented by 
Kuczynski is a
generalization of the MR (minimum residual) algorithm 
used to compute
approximations to $A^{-1}b$. There one seeks vectors $x$ 
in $R^n$ s.t.
$\|Ax-b\|<\varepsilon \|b\|$. A well-chosen subspace may 
be used to generate
approximate solutions at low cost. It is advisable to 
ensure that the
right-hand side is in the chosen subspace, and this 
consideration leads one to
choose the subspace $K^j$. In this context $b$ is part of 
the data $(A,b)$ and
is not at our disposal. The study of bad $b$'s is 
relevant to a study of the
complexity of Krylov space methods for {\it linear 
equations\/}. However, it has
been appreciated from the beginning that for reasonable 
$\varepsilon $ and
unfortunate $b$ then $n$ steps will be required unless 
$A$ is close to the
identity matrix (see [Ka, 1966; Thr. 4.3; and Me, 1963].
\par To burden the Lanczos algorithm (or GMR) with {\it 
unnecessarily\/} bad
starting vectors for the eigenvalue problem is like 
studying the performance of
Olympic athletes only when they suffer from some rare 
affliction like poison
ivy.
\subheading{{\rm 5.3.}\quad Redefining the Lanczos 
algorithm} The new algorithm
GMR is contrasted with the well-known Lanczos algorithm. 
Here is Kuczynski's
definition of the Lanczos algorithm, from p.\ 142. The 
subspace $A_j$ is our
Krylov subspace $K^j$.
\par Perform the following steps:
\roster
\item"(1)" Find an orthonormal basis $q_1,q_2,\dotsc, 
q_j$ of the subspace
$A_j$; let $Q_j=(q_1,\dotsc, q_j)$ be the $n\times j$ 
matrix.
\item"(2)" Form the $j\times j$ matrix $H_j=Q^t_jAQ_j$; 
compute eigenpairs of
$H_j$; $H_jg_i=h_ig_i$, $(g_i,g_m)=\delta_{im}$, $i, 
m=1,\dotsc, j$.
\item"(3)" Compute the Ritz vectors $z_i=Q_jg_i$ and the 
residual
$$r^L_j=\operatorname{min}_{1\le i\le 
j}\|Az_i-\theta_iz_i\|\quad\text{for }1\le i\le j.$$
\item"(4)" Define $Z_j=\{z_i,\theta_i)$, $i=1,2,\dotsc, j\:
\|Az_i-z_i\theta_i\|=r^L_j\}$
\endroster
The $j$th step of the $L$ algorithm is defined by
$$\Phi^L_j(N_j(A,b))=(x_k,\rho_k),$$
where $(x_k,\rho_k)$ is an arbitrary element from $Z_j$.
\par The trouble is that steps 3 and 4 have been changed 
from the usual ones to
conform with the idiosyncratic goal discussed in \S5.1. 
However, no mention of
this fact is made.
\par Here is the conventional description wherein it is 
supposed that $p$
eigenvalues are to be approximated. It is from \cite{Pa, 
1980, p. 214}.
\roster
\item"(1)" Is the same as above.
\item"(2)" Form the $j\times j$ matrix $H_j=Q^t_jAQ_j$; 
compute the $p(\le j)$
eigenpairs of $H_j$ that are of interest, say
$$H_jg_i=g_i\theta_i,\qquad i=1,\dotsc, p$$
The $\theta_i$ are Ritz values, 
$\theta_1<\theta_2<\dotsb< \theta_j$. Equality
is not possible.
\item"(3)" If desired, compute the $p$ Ritz vectors 
$z_i=Q_jg_i$, $i=1,\dotsc,
p$. The full set $\{(\theta_i,z_i), i=1,\dotsc, j\}$ is 
the best set of $j$
approximations to eigenpairs of $A$ that can be derived 
from $A_j$ alone.
\item"(4)" Residual error bounds.
Form the $p$ residual vectors $r_i=Az_i-z_i\theta_i$. 
Each interval
$[\theta_i-\|r_i\|, \theta_i+\|r_i\|]$ contains an 
eigenvalue of $A$. If some
intervals overlap then a bit more work is required to 
guarantee approximations
to $p$ eigenvalues. See [Pa, 1980, \S11-5]. 
\item"(5)" If satisfies then stop.
\endroster
\par In the context of Kuczynski's investigations his 
modification is entirely
reasonable; i.e., he selects at each step one Ritz pair 
with a minimal residual
norm. However, it is most misleading to those not 
familiar with the Lanczos
algorithm to suggest that its purpose is simply to 
product this Ritz pair. In
fact, as indicated above, the Lanczos algorithm produces 
an approximation,
albeit crude, to the whole spectrum, namely 
$\theta_1,\dotsc, \theta_j$ and the
user is free to select from this set to suit the specific 
goal. Thus, 
to approximate the right-most eigenvalue, one 
concentrates on $\theta_j$ and
continues until its error bound $\|r_j\|$ is 
satisfactory. In practice, more
refined error bounds can be made but that is not germane 
here (see [Pa \& No,
1985]).
\par It would have been preferable to state the Lanczos 
algorithm
conventionally and then specify the modifications 
appropriate for the purpose
in hand. This action would make clear that the Lanczos 
algorithm is not trying
to minimize one residual norm. That is why it is inferior 
to GMR for that
purpose.
\par It is worth pointing out here, that in the model of 
arithmetic used in
these studies, the cost of all the Rayleigh-Ritz 
approximations and of finding
the minimal residual norm over the subspace $K^j$ is 
taken as nil. It might
occur to the reader that in this context it would cost no 
more per step to
compute {\it all\/} the Rayleigh-Ritz approximations and 
use whatever
approximations one desires. Thus setting up GMR and 
Lanczos as competing
algorithms is artificial. Moreover, in practice, it is 
much more expensive to
compute the minimal residual than to compute the 
Rayleigh-Ritz residuals.
Kuczynski has devoted a whole report to the task.
\subheading{{\rm 5.4.}\quad Theoretical results} From p.\ 
138 of [Ku, 1986]. 
We are ready to formulate the main theory of the paper.
\proclaim{Theorem 3.1} If $F$ is unitarily 
\RM(orthogonally\RM) invariant,
then the GMR is almost
 strongly optimal in $F$, i.e., $k(\Phi^{gmr}, A, b)=
\operatorname{min}_\Phi k(\Phi,A,b)\,+a$, for any 
$(A,b)\in F\times S_n$, where
$a\in \{0,1,2\}$.
\endproclaim
\par Here $k(\Phi,A,b)$ is the minimal number of steps 
$j$ required to
guarantee an $\varepsilon $-residual with algorithm 
$\Phi$ over all matrices
$\widetilde{A}$ that are indistinguishable from $A$ with 
the given information
$N=N_j=[b,Ab,\dotsc, A^jb]$. The algorithm $\Phi$ returns 
a pair $\rho$, $x$
whose residual norm $\|\widetilde{A}x-x\rho\|$ is to be 
compared with
$\varepsilon $.
\par Recall that with information $N_j$, the GMR 
algorithm, by definition,
picks out a unit vector $x\in K^j$ and a $\rho\in C$ that 
produce the minimal
residual norm.
\par How could any other algorithm possibly do better? 
Well, there might be
special circumstances in which one could deduce a 
suitable additional component
of that last vector $A^jb$ that is not used by GMR in 
forming $x$, although
$A^jb$ is used in calculating the coefficients of GMR's 
approximation from
$K^j$. The proof studies this possibility and concludes 
that GMR would make up
any discrepancy in at most two more steps.
\par The argument is very nice. In many important cases, 
when $A$ is Hermitian
for example, then the constant $a$ in Theorem 3.1 is 
actually 0.
\par There are other clever results. Theorem 4.2 shows 
that for symmetric 
matrices the residual norm of GMR must be strictly 
decreasing at lest at every
other step. Theorem 5.1 yields a beautiful but esoteric 
fact about Krylov
subspaces generated by Hermitian $A$. For the worst 
starting vector there 
is a unit vector $v$ in $K^j$ such that
$$\frac{\|A\|}{2j}\le \|Av-v\rho\|\le\frac{\|A\|}j,$$
for $j<n$.
\par As indicated above, these nice results from 
approximation theory do not
add up to a case for replacing Rayleigh-Ritz with some 
rival algorithm.
We sometimes abbreviate Rayleigh-Ritz by R-R.
\subheading{{\rm 5.5}.\quad Numerical examples} All the 
numerical results
reported in [Ku, 1986] concern symmetric tridiagonal 
matrices with starting
vector $e_1$ (the first column of the identity matrix 
$I)$. This starting
vector ensures that the Lanczos algorithm reproduces the 
original matrix. At
this point we should recall that the original goal of the 
Lanczos algorithm was
to reduce a symmetric matrix to tridiagonal form. So, the 
numerical 
results to be seen below do not relate to the Lanczos 
recurrence itself but
merely indicate alternative rules for stopping. With the 
goal of
tridiagonalization, the algorithm always stopped at step 
$n$. Later, it was 
realized that excellent approximations to a few 
eigenvectors were usually
obtained early in the process. There is no single 
stopping criterion for
current Lanczos algorithms; termination depends on what 
the user wants (see [Cu
\& Wi, 1985; Pa, 1980; and Go \& Va, 1984].
\par Kuczynski provides the Lanczos algorithm with a 
stopping criterion to suit
his purposes, but his algorithm GMR could have been 
called (with more justice)
the Lanczos algorithm with a new stopping criterion. It 
uses the minimum
residual in the whole Krylov subspace instead of the 
usual (cheap)
Rayleigh-Ritz approximations. So, the numerical results 
simply indicate the
effect of these different termination criteria.
\par The first batch of results concern tridiagonals with 
nonzero elements
chosen at random from $[-\tfrac13, \tfrac13]$. The most 
striking feature is
that the GMR residual and the smallest Rayleigh-Ritz 
residual slip below the
given $\varepsilon $ at the same step in the vast 
majority of
cases, particularly for $\varepsilon <10^{-3}$. In Table 
8.1, with $\varepsilon
=10^{-6}$, the step was the same in 18 out of 20 cases. 
In the other two, the
Lanczos algorithm (i.e. the minimal $R\text{-}R$ norm) 
took one more step (17 as
against 16).
\par Some weight is given to the fact that the smallest 
$R\text{-}R$ residual norm is
rarely monotone decreasing from one step to another 
whereas GMR does enjoy this
property. However, if the approximate eigenvalue 
associated with the minimum
residual happens to change position in the spectrum from 
step to step, then
this monotonicity of GMR is not associated with the 
convergence to a specific
eigenvalue of the original matrix. No indication is given 
in the results of how
the approximate eigenvalue implicitly chosen by GMR jumps 
around the spectrum.
In practice, the interesting thing to know is {\it how 
many\/} Ritz
values have {\it converged\/}, and to what accuracy, when 
the algorithm is
terminated. Unfortunately, this information is excluded 
from the GMR viewpoint 
and is not reported.
\par The next results, Examples 8.1 and 8.2 in [Ku, 
1986], exhibit the
dramatic {\it failure\/} of the Lanczos algorithm. On a 
tridiagonal matrix of
order 201 and norm near 1 the minimal $R\text{-}R$ 
residual remained at its initial
value 0.035 for all steps except the last (at which it 
must be 0). In contrast,
the GMR residual declined slowly from the intial 0.035 to 
0.0039 at step 200.
If $\varepsilon =0.034$, then GMR takes 2 steps while 
Lanczos takes 201!
However, with $\varepsilon \le 10^{-3}$, both algorithms 
need 201 steps. We
repeat, once again, that GMR will not know which 
eigenvalue it has
approximated.
\par Unfortunately, no attempt is made to put this 
example in context. It
illustrates the phenomenon explored in some detail in 
[Sc, 1979], namely that
for every symmetric matrix with distinct eigenvalues 
there is a set (with
positive Lebesgue measure on the sphere) of starting 
vectors such that no
Rayleigh-Ritz approximation is any good until the last 
step. We must repeat
that the Lanczos algorithm is not obliged to use a poor 
initial vector. We ran
{\it our\/} Lanczos code on this matrix. Our code starts 
with a normalized
version of $Ar$, where $A$ is the given matrix (or 
operator) and $r$'s elements
are chosen at random from a uniform random distribution. 
The reason for
starting with $Ar$ is compelling when an operator $A$
 has unwanted infinite
eigenvalues. The results are given in the
 Table 1 (see p.\ 24).
\midinsert
\toptablecaption{{\smc Table 1}. Convergence of Ritz 
pairs on $T_{201}$.}
\noindent$$
\table\format\c&\rulev\quad\c&\rulev\quad\c\\
\text{STEP}&\varepsilon&\colheading{Number of\\
good Ritz values}\halfspace\\
\ruleh
20&10^{-2}&1\\
30&10^{-3}&1\\
40&10^{-3}&4\\
50&10^{-3}&5\\
60&10^{-4},\ 10^{-6}&7,4\\
70&10^{-5},\ 10^{-7}&10,5\\
80&10^{-5},\ 10^{-7}&14,10\\
90&10^{-5},\ 10^{-7}&18,14\\
100&10^{-5},\ 10^{-7}&23,18\\
110&10^{-5},\ 10^{-7}&28,24\\
120&10^{-5},\ 10^{-7}&35,30\\
130&10^{-5},\ 10^{-7}&44,37\\
140&10^{-5},\ 10^{-7}&52,44\\
150&10^{-5},\ 10^{-7}&60,54\\
160&10^{-5},\ 10^{-7}&70,61\\
170&10^{-5},\ 10^{-7}&83,75\\
180&10^{-5},\ 10^{-7}&99,89\\
190&10^{-5},\ 10^{-7}&118,109
\endtable$$
\endinsert
\par The accepted eigenvalues (104 of them at step 190) 
agreed with those
computed by EISPACK to all of the fifteen decimals 
printed out. The efficiency is
not at all bad considering that this is a difficult 
eigenvalue distribution for
Krylov space methods.
\par Example 8.2, a tridiagonal of order 501 with null 
diagonal and monotonely
increasing off diagonal elements, caused the minimal 
$R\text{-}R$ residual norm to
{\it increase\/} from 0.001 initially in 0.011 at steps 
499 and 500. In
contrast GMR residual norms declined to 0.00036 at steps 
499 and 500. Thus,
with $\varepsilon =.00099$, GMR terminates at step 2 
whereas Lanczos terminates
at step 501! However, with $\varepsilon \le 10^{-4}$, 
both take 501 steps.
\par As with Example 8.1 $e_1$ is a bad staring vector 
yielding a poor Krylov
subspace. We ran our Lanczos program and found the 
results given in Table 2.
\midinsert
\toptablecaption{{\smc Table 2}. Convergence of Ritz 
pairs on $T_{501}$.}
\noindent$$
\table\format\c\rulev&\quad\c\rulev&\quad\c\\
\text{STEP}&\varepsilon &\colheading{Number of\\
good Ritz values}\halfspace\\
\ruleh
20&10^{-2}&1\\
30&10^{-3}&1\\
40&10^{-3}&2\\
50&10^{-4}&3\\
60&10^{-4},\ 10^{-6}&5,2\\
70&10^{-5},\ 10^{-7}&7,3\\
80&10^{-5},\ 10^{-7}&9,6\\
90&10^{-5},\ 10^{-7}&13,10\\
100&10^{-5},\ 10^{-7}&16,13\\
110&10^{-5},\ 10^{-7}&20,15\\
120&10^{-5},\ 10^{-7}&23,20\\
130&10^{-5},\ 10^{-7}&29,23\\
140&10^{-5},\ 10^{-7}&34,29\\
150&10^{-5},\ 10^{-7}&39,34
\endtable$$
\endinsert
\par \phantom{The\ }We quote the final paragraph of the 
article.
\ext From all the tests we have performed we conclude 
that the GMR
algorithm is
essentially superior to the Lanczos Algorithm on matrices 
with constant or
increasing codiagonal elements. For random matrices or 
matrices with decreasing
codiagonal elements, both algorithms produce nearly the 
same residuals.
\endext
\par The revealing word here is ``codiagonal.'' The 
author has worked
exclusively with tridiagonal matrices and has forgotten 
that the goal of the
Lanczos recurrence is to produce a tridiagonal matrix! 
Given such a matrix one
has {\it no need of either Lanczos or \/} GMR. As our 
results indicate, a
random starting vector permits the Lanczos algorithm to 
perform satisfactorily
even on such craftily-designed matrices. The quotation 
reveals just how far a
mathematical theory can stray from relevance.
\subheading{{\rm 5.6.}\quad Summary} Here is an attempt 
to formulate the
numerical analyst's version of Complexity Theory for 
Krylov subspaces and
eigenvalues.
\ext
For each symmetric $n\times n$ matrix there are initial 
vectors that yield an
eigenvalue in one step, and initial vectors that yield an 
eigenvalue only at
the $n$th step. The nontrivial result contained in the 
Kaniel-Paige-Saad error
bounds (see [Pa, 1980, Chap.\ 12]) is that with most 
starting vectors the 
extreme eigenvalues can be found in a modest number of 
steps that depends on
the distribution of the spectrum and is nearly 
independent of $n$.
\endext
\par We summarize our criticism of [Ku, 1986] but wish to 
note that the paper
is essentially the author's Ph.D. dissertation, and it 
would not be charitable
to hold him responsible for the vagaries to which IBCT is 
subject.
\varroster
\item""Section 1 exposes a serious flaw in the model, 
namely
the goal.
\item""Section 2 exposes a subtle way in which features 
of a method are pushed
into the problem statement; the starting vector.
\item""Section 3 shows how standard terms can be 
redefined; the Lanczos 
algorithm.
\item""Section 4 contains some clever and interesting 
results on the
approximating power of Krylov subspaces.
\item""Section 5 shows how very misleading numerical 
results can be in the
absence of proper interpretation.
\endroster
\heading Acknowledgment\endheading
\par The author has been assisted by many friends in 
composing this essay but
would like to single out Stan Eisenstat, G. W. Stewart, 
and W. Kahan for their
helpful and detailed comments.


\Refs
\widestnumber\key{O'L, Ste \& Va, 1979}
\ref\key Ba, 1987 \by I. Babuska \paper Information-based 
numerical
practice \jour J. Complexity \vol 3 \yr 1987 \pages 
331--346\endref
\ref\key Bl, Shu, \& Sm, 1989 \by L. Blum, M. Shub, and 
S. Smale
\paper On a theory of computation and complexity over the 
real numbers$\!;$
NP completeness, recursive functions, and Turing machines 
\jour Bull.
Amer. Math. Soc. (N.S.) \vol 21 \yr 1989 \pages 
1--46\endref
\ref\key Bo \& Mu, 1975 \by A. Borodin and I. Munro \book 
Computational
complexity of algebraic and numeric problems \publ Amer. 
Elsevier,
New York, 1975\endref
\ref\key Cu \& Wi, 1985 \by J. K. Cullum and R. A. 
Willoughby \book
Lanczos algorithms for large symmetric eigenvalue 
computations,
Vol. {\rm I:} Theory \bookinfo Progr. Comput. Sci., 
Birkhauser,
Basel, 1985\endref
\ref\key Da, 1967 \by J. W. Daniel \paper The conjugate 
gradient
method for linear and nonlinear operator equations \jour 
SIAM
J. Numer. Anal. \vol 4 \yr 1967 \pages 10--26\endref
\ref\key Go \& Va, 1984 \by G. H. Golub and C. V. Van 
Loan \book
Advanced matrix computations \bookinfo Johns Hopkins 
Univ. Press,
Maryland, 1984\endref
\ref\key Je, 1977 \by A. Jennings \book Matrix 
computation for
engineers and scientists \bookinfo John Wiley \& Sons, 
Chichester, 1977
\endref
\ref\key Ka, 1966 \by S. Kaniel \paper Estimates for some 
computational
techniques in linear algebra \jour Math. Comp. \vol 20 
\yr 1966
\pages 369--378\endref
\ref\key Ku, 1986 \by J. Kuczynski \paper A generalized 
minimal residual
algorithm for finding an eigenpair of a large matrix 
\jour J.
Complexity \vol 2 \yr 1986 \pages 131--162\endref
\ref\key Ku, 1985 \bysame \book Implementation of the GMR 
algorithm
for large symmetric eigenproblems \bookinfo Tech. Report, 
Comp. Sci.
Dept., Columbia University, New York, 1985\endref
\ref\key Ku \& Wo, 1992 \by J. Kuczynski and H. 
Wo\'zniakowski
\paper Average case complexity of the symmetric Lanczos 
algorithm
\jour SIAM J. Matrix Anal. Appl., 1992 (to appear)\endref
\ref\key La, 1952 \by C. Lanczos \paper Solution of 
systems of linear
equations by minimized iterations \jour J. Res. Nat. Bur. 
Standards
\vol 49 \yr 1952 \pages 33--51\endref
\ref\key Me, 1963 \by G. Meinardus \paper \"Uber eine 
Verallgemeinerung einer
Ungleichung von L. V. Kantorowitsch \jour Numer. Math. 
\vol 5
\yr 1963 \pages 14--23\endref
\ref\key O'L, Ste \& Va, 1979 \by D. O'Leary, G. W. 
Stewart and J. S.
Vandergraft \paper Estimating the largest eigenvalue of a 
positive
definite matrix \jour Math. Comp. \vol 33 \yr 1979 \pages 
1289--1292
\endref
\ref\key Pa, 1980 \by B. N. Parlett \paper The symmetric 
eigenvalue
problem \jour Prentice-Hall, New Jersey, 1980\endref
\ref\key Pa, 1982 \bysame \book Two monitoring schemes 
for the Lanczos
algorithm \bookinfo Computing Methods in Applied Sciences 
and Engineering
(V. R. Glowinski and J. L. Lions, eds.), North-Holland 
Press, 1982\endref
\ref\key Pa, 1984 \bysame \paper The software scene in 
the extraction
of eigenvalues from sparse matrices \jour SIAM J. Sci. 
Statist. Comput.
\vol 5 \yr 1984 \pages 590--603\endref
\ref\key Pa \& No, 1985 \by B. N. Parlett and B. 
Nour-Omid \paper
The use of refined error bounds when updating eigenvalues 
of tridiagonals
\jour Linear Algebra Appl. \vol 68 \yr 1985 \pages 
179--219\endref
\ref\key Pa, Si \& Str, 1982 \by B. N. Parlett, H. A. 
Simon, and L. M.
Stringer \paper On estimating the largest eigenvalue with 
the Lanczos
algorithm \jour Math. Comp. \vol 38 \yr 1982 \pages 
153--165\endref
\ref\key Pac, 1986 \by E. Packel \paper {\rm Review of} A 
general theory of
optimal algorithms \jour by J. F. Traub and H. 
Wo\'zniakowski (Academic
Press, New York, 1980), SIAM Rev. \vol 28 \yr 1986 \pages 
435--437\endref
\ref\key Pac \& Wo, 1987 \by E. Packel and H. 
Wo\'zniakowski \paper
Recent developments on information-based complexity \jour 
Bull. Amer.
Math. Soc. (N.S.) \vol 17 \yr 1987 \pages 9--26\endref
\ref\key Ra, 1972 \by M. O. Rabin \paper Solving linear 
equations by means
of scalar products \inbook Complexity of Computer 
Computations
(R. E. Miller and J. W. Thatcher, eds.), Plenum Press, 
New York, 1972, pp.
11--20\endref
\ref\key Sa, 1980 \by Y. Saad \paper On the rates of 
convergence of the
Lanczos and the block-Lanczos methods \jour SIAM J. 
Numer. Anal.
\vol 17 \yr 1980 \pages 687--706\endref
\ref\key Sc, 1979 \by D. S. Scott \paper How to make the 
Lanczos algorithm
converge slowly \jour Math. Comp. \vol 33 \yr 1979 \pages 
239--247\endref
\ref\key Sch \& Str, 1971 \by A. Schoenhage and V. 
Strassen \paper
Schnelle Multiplikation Grosser Zahlen \jour Commuting 
\vol 7
\yr 1971 \pages 281--292\endref
\ref\key Sh, 1977 \by I. Shavitt \paper The method of 
configuration
interaction \jour Modern Theoretical Chemistry, vol. 3 
(H. P.
Schaefer, ed.), Plenum Press, 1977\endref
\ref\key Shu, 1987 \by M. Shub \paper {\rm Review of} 
Information,
uncertainty, complexity \jour by J. F. Traub, H. 
Wo\'zniakowski et al.
(Addison-Wesley, Reading, MA 1983), SIAM Rev. \vol 29 \yr 
1987
\pages 495--496\endref
\ref\key Sti, 1958 \by E. Stiefel \paper Kernel 
polynomials in linear
algebra and their numerical applications \jour NBS Appl. 
Math.
\vol 43 \yr 1958 \pages 1--22\endref
\ref\key Str, 1969 \by V. Strassen \paper Gaussian 
elimination is not
optimal \jour Numer. Math. \vol 13 \yr 1969 \pages 
354--356\endref
\ref\key Tr \& Wo, 1980 \by J. F. Traub and H. 
Wo\'zniakowski
\paper A general theory of optimal algorithms \jour 
Academic Press, New
York, 1980\endref
\ref\key Tr \& Wo, 1984 \bysame \paper On the optimal 
solution of large
linear systems \jour J. Assoc. Comput. Mach., \vol 31 \yr 
1984
\pages 545--559\endref
\ref\key Tr \& Wo, 1988 \bysame \paper Information-based 
complexity
\jour Academic Press, New York, 1988\endref
\ref\key Tr, Wo \& Wa, 1983 \by J. F. Traub, H. 
Wo\'zniakowski, and G.
Wasilkowski \paper Information, uncertainty, complexity 
\jour Addison-Wesley,
Reading, 1983\endref
\ref\key Wi, 1970 \by S. Winograd \paper On the number of 
multiplication
necessary to compute certain functions \jour Comm. Pure 
Appl. Math.
\vol 23 \yr 1970 \pages 165--179\endref
\ref\key Wi, 1980 \bysame \paper Arithmetic complexity of 
computations
\jour CBMS-NSF Regional Conf. Ser. In Appl. Math., vol. 
33, SIAM,
Philadelphia, 1980\endref
\ref\key Wo, 1985 \by H. Wo\'zniakowski \paper A survey 
of information-based
complexity \jour J. Complexity \vol 1 \yr 1985 \pages 
11--44\endref



\endRefs
\enddocument